\documentclass[10pt]{amsart}
\usepackage{amsfonts}
\usepackage{ifthen}
\usepackage{amsthm}
\usepackage{amsmath}
\usepackage{graphicx}
\usepackage{amscd,amssymb,amsthm}

\newcounter{minutes}
\divide\time by 60
\newcounter{hours}
\multiply\time by 60 \addtocounter{minutes}{-\time}

\newtheorem{theorem}{Theorem}
\newtheorem{lemma}{Lemma}


\title[]{Close-to-convexity of normalized Dini functions}

\author[\'A. Baricz]{\'Arp\'ad Baricz$^{\bigstar}$}
\address{Department of Economics, Babe\c{s}-Bolyai University, Cluj-Napoca 400591, Romania}
\address{Institute of Applied Mathematics, John von Neumann Faculty of Informatics, \'Obuda University, 1034 Budapest, Hungary}
\email{bariczocsi@yahoo.com}

\author[E. Deniz]{Erhan Deniz}
\address{Department of Mathematics, Faculty of Science and Letters, Kafkas
University, Kars 36100, Turkey.} \email{edeniz36@gmail.com}

\author[N. Yagmur]{Nihat Yagmur}
\address{Department of Mathematics, Erzincan University, Erzincan 24000,
Turkey} \email{nhtyagmur@gmail.com}

\thanks{$^{\bigstar}$The research of \'A. Baricz was supported by a research grant of the Romanian National Authority for Scientific Research, CNCS-UEFISCDI, project number PN-II-RU-TE-2012-3-0190/2014.}

\begin{document}

\def\thefootnote{}
\footnotetext{ \texttt{File:~\jobname .tex,
          printed: \number\year-0\number\month-\number\day,
          \thehours.\ifnum\theminutes<10{0}\fi\theminutes}
} \makeatletter\def\thefootnote{\@arabic\c@footnote}\makeatother

\maketitle

\begin{center}
{\em Dedicated to Aysima, Bor\'oka and Kopp\'any}
\end{center}

\begin{abstract}
In this paper necessary and sufficient conditions are deduced for the close-to-convexity of some special combinations of Bessel functions of the first kind and their derivatives by using a result of Shah and Trimble about transcendental entire functions with univalent derivatives and some newly discovered Mittag-Leffler expansions for Bessel functions of the first kind.
\end{abstract}

\section{\bf Introduction and the Main Results}

Special functions, like Bessel functions of the first kind play an important role in pure and applied mathematics. Geometric properties, like univalence, starlikeness, spirallikeness and convexity were studied in the sixties by Brown \cite{brown, brown2,brown3}, and also by Kreyszig and Todd \cite{todd}. Some other geometric properties of Bessel functions of the first kind were studied later in the papers \cite{mathematica,publ,lecture,bsk,samy,basz,szasz,szasz2}. Recently, in \cite{barsza} the close-to-convexity of the derivatives of Bessel functions was considered. In this paper we make a contribution to the subject by deducing necessary and sufficient conditions for the close-to-convexity of some special combinations of Bessel functions of the first kind and their derivatives. In order to prove our main results we use a result of Shah and Trimble \cite[Theorem 2]{st} about transcendental entire functions with univalent derivatives and some newly discovered Mittag-Leffler expansions for Bessel functions of the first kind. The next result of Shah and Trimble \cite[Theorem 2]{st} is the cornerstone of this paper.

\begin{lemma}\label{lem1}
Let $\mathbb{D}=\{z\in\mathbb{C}:|z|<1\}$ be the open unit disk and $f:\mathbb{D}\to\mathbb{C}$ be a transcendental entire function of the form
$$f(z)=z\prod_{n\geq 1}\left(1-\frac{z}{z_n}\right),$$
where all $z_n$ have the same argument and satisfy $|z_n|>1.$ If $f$ is univalent in $\mathbb{D},$ then
$$\sum_{n\geq1}\frac{1}{|z_n|-1}\leq 1.$$
In fact the above inequality holds if and only if $f$ is starlike in $\mathbb{D}$ and all of its derivatives are
close-to-convex there.
\end{lemma}

Now, consider the function $f_{\nu}:\mathbb{D}\to\mathbb{C},$ defined by
$$f_{\nu}(z)=2^{\nu}\Gamma(\nu+1)z^{1-\frac{\nu}{2}}J_{\nu}(\sqrt{z})=\sum_{n\geq0}\frac{(-1)^n\Gamma(\nu+1)z^{n+1}}{4^nn!\Gamma(\nu+n+1)},$$
where $J_{\nu}$ stands for the Bessel function of the first kind (see \cite[p. 217]{nist}). Very recently by using a result of Shah and Trimble \cite[Theorem 2]{st} in \cite{barsza} the authors proved that the function $f_{\nu}$ and all of its derivatives are convex in $\mathbb{D}$ if and only if $\nu\geq\nu_{\star},$ where $\nu_{\star}\simeq-0.1438\dots$ is the unique root of the equation $$3J_{\nu}(1)+2(\nu-2)J_{\nu+1}(1)=0$$ on $(-1,\infty).$ We note that in view of the Alexander's duality theorem the first part of the above result is equivalent to the fact that the function
\begin{align*}z\mapsto q_{\nu}(z)=zf_{\nu}'(z)&=2^{\nu-1}\Gamma(\nu+1)z^{1-\frac{\nu}{2}}\left((2-\nu)J_{\nu}(\sqrt{z})+\sqrt{z}J_{\nu}'(\sqrt{z})\right)\\
&=\sum_{n\geq0}\frac{(-1)^n(n+1)\Gamma(\nu+1)z^{n+1}}{4^nn!\Gamma(\nu+n+1)}\end{align*}
is starlike in $\mathbb{D}$ if and only if $\nu\geq\nu_{\star}.$ In this paper we would like to point out that actually we have the following stronger result.

\begin{theorem}\label{th0}
The function $q_{\nu}$ is starlike and all of its derivatives are close-to-convex (and hence univalent) in $\mathbb{D}$ if and only if $\nu\geq\nu_{\star},$ where $\nu_{\star}\simeq-0.1438\dots$ is the unique root of the transcendent equation $3J_{\nu}(1)+2(\nu-2)J_{\nu+1}(1)=0$ on $(-1,\infty).$
\end{theorem}

Moreover, in this paper we are interested on the normalized Dini function $r_{\nu}:\mathbb{D}\to\mathbb{C},$ which is another special combination of Bessel functions of the first kind and is defined as
\begin{align*}r_{\nu}(z)&=2^{\nu}\Gamma(\nu+1)z^{1-\frac{\nu}{2}}\left((1-\nu)J_{\nu}(\sqrt{z})+\sqrt{z}J_{\nu}'(\sqrt{z})\right)\\&
=\sum_{n\geq0}\frac{(-1)^n(2n+1)\Gamma(\nu+1)z^{n+1}}{4^nn!\Gamma(\nu+n+1)}.\end{align*}

By using Lemma \ref{lem1} and the idea of the proof of Theorem \ref{th0} our aim is to present the following interesting sharp result. We note that some similar results were proved for Bessel functions of the first kind in \cite{bsk,basz,barsza,szasz}, but by using different approaches.

\begin{theorem}\label{th1}
The function $r_{\nu}$ is starlike and all of its derivatives are close-to-convex (and hence univalent) in $\mathbb{D}$ if and only if $\nu\geq\nu^{\star},$ where $\nu^{\star}\simeq0.3062\dots$ is the unique root of the transcendent equation $J_{\nu}(1)-(3-2\nu)J_{\nu+1}(1)=0$ on $(0,\infty).$
\end{theorem}

It is important to mention here that the first part of the above result is quite similar to the following sharp result (see \cite[Theorem 1.5]{basz}): the function $z\mapsto 2^{\nu}\Gamma(\nu+1)z^{1-\nu}J_{\nu}(z)$ is convex in $\mathbb{D}$ if and only if $\nu\geq1.$ Note that according to the Alexander's duality theorem this result is equivalent to the following: the function $$z\mapsto 2^{\nu}\Gamma(\nu+1)z^{1-\nu}\left((1-\nu)J_{\nu}(z)+zJ_{\nu}'(z)\right)=\frac{r_{\nu}(z^2)}{z}$$ is starlike in $\mathbb{D}$ if and only if $\nu\geq1.$

Now, let us consider the function $w_{a,\nu}:\mathbb{D}\to \mathbb{C},$ defined by
$$w_{a,\nu}(z)=\frac{2^{\nu}}{a}\Gamma(\nu+1)z^{1-\frac{\nu}{2}}\left((a-\nu)J_{\nu}(\sqrt{z})+\sqrt{z}J_{\nu}'(\sqrt{z})\right)=
\sum_{n\geq0}\frac{(-1)^n(2n+a)\Gamma(\nu+1)z^{n+1}}{a\cdot 4^nn!\Gamma(n+\nu+1)}.$$
The following sharp result is a common generalization of Theorems \ref{th0} and \ref{th1}.

\begin{theorem}\label{thg}
Let $\nu>-\frac{3}{4}$ and $a\geq \frac{2}{4\nu+3}.$ The function $w_{a,\nu}$ is starlike and all of its derivatives are close-to-convex (and hence univalent) in $\mathbb{D}$ if and only if $\nu\geq\nu_a,$ where $\nu_a$ is the unique root of the transcendent equation $(2a-1)J_{\nu}(1)-(a-2\nu+2)J_{\nu+1}(1)=0$ on $\left(-\frac{3}{4},\infty\right).$
\end{theorem}

Finally, we mention that in particular Theorem \ref{th0} and Theorem \ref{th1} yield that $$z\mapsto q_{\frac{1}{2}}(z)=\frac{3}{2}\sqrt{z}\left(\sin\sqrt{z}+\sqrt{z}\cos\sqrt{z}\right),$$
$$z\mapsto q_{\frac{3}{2}}(z)=\frac{3}{2\sqrt{z}}\left(\sqrt{z}\cos\sqrt{z}+(z-1)\sin\sqrt{z}\right),$$
$$z\mapsto r_{\frac{1}{2}}(z)=z\cos\sqrt{z}$$ and $$z\mapsto r_{\frac{3}{2}}(z)=3\cos\sqrt{z}-\frac{3(z-2)\sin\sqrt{z}}{2\sqrt{z}}$$ are starlike in $\mathbb{D}$ and all of their derivatives are close-to-convex (and hence univalent) there. Moreover, by using the above result we obtain that $$z\mapsto \frac{r_{\frac{3}{2}}(z^2)}{z}=\frac{3\cos z}{z}-\frac{3(z^2-2)\sin z}{2z^2}$$ is starlike in $\mathbb{D}.$ Here we used that \cite[p. 228]{nist}
$$J_{\frac{1}{2}}(z)=\sqrt{\frac{2}{\pi z}}\sin z\ \ \ \mbox{and}\ \ \ J_{\frac{3}{2}}(z)=\sqrt{\frac{2}{\pi z}}\left(\frac{\sin z}{z}-\cos z\right).$$

\section{\bf Proofs of the Main Results}
\setcounter{equation}{0}

In this section our aim is to present the proof of our main results.

\begin{proof}[\bf Proof of Theorem \ref{th0}]
Let us denote the $n$th positive zero of the Dini function $z\mapsto (2-\nu)J_{\nu}(z)+zJ_{\nu}'(z)$ by $\beta_{\nu,n}.$ We know that \cite[Lemma 2.5]{basz} if $\nu>-1,$ then we have the Mittag-Leffler expansion
$$\frac{zf_{\nu}''(z)}{f_{\nu}'(z)}=-\sum_{n\geq 1}\frac{z}{\beta_{\nu,n}^2-z},$$
which in turn implies that
$$f_{\nu}'(z)=e^{c}\prod_{n\geq1}\left(1-\frac{z}{\beta_{\nu,n}^2}\right).$$
On the other hand, we know that
$$f_{\nu}'(z)=\sum_{n\geq0}\frac{(-1)^n(n+1)\Gamma(\nu+1)z^{n}}{4^nn!\Gamma(\nu+n+1)}$$
and then we have $f_{\nu}'(0)=1.$ This implies that $c=0.$ Thus for the function $q_{\nu}$ we can apply Lemma \ref{lem1}, and we just need to check the corresponding conditions on the zeros $\beta_{\nu,n}.$ We know that (see the proof of \cite[Theorem 1]{barsza}) $$\sum_{n\geq 1}\frac{1}{\beta_{\nu,n}^2-1}\leq 1$$
if and only if $\nu\geq \nu_{\star},$ where $\nu_{\star}\simeq-0.1438\dots$ is the unique root of the equation $3J_{\nu}(1)+2(\nu-2)J_{\nu+1}(1)=0$ on $(-1,\infty).$ Thus, we need only to verify the condition $\beta_{\nu,n}>1$ for each $\nu>\nu_{\star}$ and $n\in\mathbb{N}.$ In \cite[Lemma 2.6]{basz} it was proved that $\beta_{\nu,1}>1$ for $\nu\geq0,$ however a close inspection on the proof reveals that $\beta_{\nu,1}>1$ in fact for $\nu>-\frac{1}{4}.$ This implies that indeed $\beta_{\nu,n}>1$ for each $\nu>\nu_{\star}$ and $n\in\mathbb{N},$ and thus we can apply Lemma \ref{lem1} for the function $q_{\nu}.$
\end{proof}

\begin{proof}[\bf Proof of Theorem \ref{th1}]
Let us start with the following result, which was proved quite recently in \cite[Theorem 1]{dini}: Let $\nu>-1$ and consider the Dini function $d_{\nu}:\mathbb{D}\to\mathbb{C},$ defined by $$d_{\nu}(z)=(1-\nu)J_{\nu}(z)+zJ_{\nu}'(z).$$ If $\alpha_{\nu,n}$ denotes the $n$th positive zero of the Dini function $d_{\nu},$ then the following Weierstrassian factorization is valid
$$d_{\nu}(z)=\frac{z^{\nu}}{2^{\nu}\Gamma(\nu+1)}\prod_{n\geq1}\left(1-\frac{z^2}{\alpha_{\nu,n}^2}\right),$$
where the infinite product is uniformly convergent on each compact subset of the complex plane. By using this result we obtain that
$$r_{\nu}(z)=2^{\nu}\Gamma(\nu+1)z^{1-\frac{\nu}{2}}d_{\nu}(\sqrt{z})=z\prod_{n\geq 1}\left(1-\frac{z}{\alpha_{\nu,n}^2}\right).$$
On the other hand, by using the Mittag-Leffler expansion \cite[Lemma 2.4]{basz} (see also \cite[Theorem 3]{dini})
$$\frac{d_{\nu}'(z)}{d_{\nu}(z)}-\frac{\nu}{z}=\frac{zJ_{\nu+2}(z)-3J_{\nu+1}(z)}{J_{\nu}(z)-zJ_{\nu+1}(z)}=
-\sum_{n\geq1}\frac{2z}{\alpha_{\nu,n}^2-z^2}$$
we obtain that
$$\sum_{n\geq 1}\frac{1}{\alpha_{\nu,n}^2-1}=-\frac{1}{2}\cdot \frac{J_{\nu+2}(1)-3J_{\nu+1}(1)}{J_{\nu}(1)-J_{\nu+1}(1)}.$$
Now, for $\nu>-1$ let $\gamma_{\nu,n}$ be the $n$th positive root of the equation $\gamma{J}_\nu(z)+zJ_\nu'(z)=0.$ Owing to Landau \cite[p. 196]{landau} we know that if $\nu+\gamma\geq0,$ then the function $\nu\mapsto \gamma_{\nu,n}$ is strictly increasing on $(-1,\infty)$ for $n\in\mathbb{N}$ fixed. This implies that $\nu\mapsto \alpha_{\nu,n}$ is strictly increasing on $(-1,\infty)$ for $n\in\mathbb{N}$ fixed, and thus
the function $$\nu\mapsto 1-\sum_{n\geq 1}\frac{1}{\alpha_{\nu,n}^2-1}$$ is strictly increasing on $(-1,\infty).$
Consequently we have that $$\sum_{n\geq 1}\frac{1}{\alpha_{\nu,n}^2-1}\leq 1$$
if and only if $\nu\geq \nu^{\star},$ where $\nu^{\star}\simeq0.3062\dots$ is the unique root of the equation $J_{\nu}(1)-(3-2\nu)J_{\nu+1}(1)=0$ on $(0,\infty).$ Here we used the recurrence relation
$$J_{\nu}(z)+J_{\nu+2}(z)=\frac{2(\nu+1)}{z}J_{\nu+1}(z)$$
to show that the equation $$-\frac{1}{2}\left(J_{\nu+2}(1)-3J_{\nu+1}(1)\right)=J_{\nu}(1)-J_{\nu+1}(1)$$
is equivalent to $J_{\nu}(1)-(3-2\nu)J_{\nu+1}(1)=0.$ We also note that according to \cite[Lemma 2.6]{basz} for $\nu\geq0$ we have $\alpha_{\nu,1}>1$ and consequently we have $\alpha_{\nu,n}>1$ for all $n\in\mathbb{N}$ and $\nu\geq0.$ Thus, by using Lemma \ref{lem1} we conclude that the function $r_{\nu}$ is indeed starlike and all of its derivatives are close-to-convex in $\mathbb{D}$ if and only if $\nu\geq\nu^{\star}.$
\end{proof}

\begin{proof}[\bf Proof of Theorem \ref{thg}]
Let us consider the function $D_{a,\nu}:\mathbb{D}\to\mathbb{C},$ defined by
$$D_{a,\nu}(z)=(a-\nu)J_{\nu}(z)+zJ_{\nu}'(z)=\sum_{n\geq 0}\frac{(-1)^n(2n+a)z^{2n+\nu}}{2^{2n+\nu}n!\Gamma(n+\nu+1)}.$$
Since for $c>0$ we have that $\frac{\log \Gamma(n+c)}{n\log n}\to 1$ as $n\to\infty,$ it follows that
$$\lim_{n\to\infty}\frac{n\log n}{\log\left(\frac{a\cdot4^n}{(2n+a)\Gamma(\nu+1)}\right)+\log\Gamma(n+1)+\log\Gamma(n+\nu+1)}=\frac{1}{2}.$$
Thus, the growth order of the entire function
$$z\mapsto \frac{2^{\nu}}{a}\Gamma(\nu+1)z^{-\nu}D_{a,\nu}(z)=\sum_{n\geq0}\frac{(-1)^n(2n+a)\Gamma(\nu+1)z^{2n}}{a 4^nn!\Gamma(n+\nu+1)}$$
is $\frac{1}{2}$ and in view of the Hadamard theorem \cite[p. 26]{lev} we have that
$$D_{a,\nu}(z)=\frac{az^{\nu}}{2^{\nu}\Gamma(\nu+1)}\prod_{n\geq 1}\left(1-\frac{z^2}{\omega_{a,\nu,n}^2}\right),$$
where $\omega_{a,\nu,n}$ stands for the $n$th positive zero of $D_{a,\nu}.$ This in turn implies that
$$w_{a,\nu}(z)=\frac{2^{\nu}}{a}\Gamma(\nu+1)z^{1-\frac{\nu}{2}}D_{a,\nu}(\sqrt{z})=z\prod_{n\geq1}\left(1-\frac{z}{\omega_{a,\nu,n}^2}\right).$$
On the other hand, by using logarithmic differentiation in the infinite product representation of $D_{a,\nu}$ and the fact that the Bessel function of the first kind is a particular solution of the second-order Bessel differential equation, it follows the Mittag-Leffler expansion
$$\frac{D_{a,\nu}'(z)}{D_{a,\nu}(z)}-\frac{\nu}{z}=\frac{(2\nu^2-a\nu-z^2)J_{\nu}(z)+(a-2\nu)zJ_{\nu}'(z)}{(a-\nu)zJ_{\nu}(z)+z^2J_{\nu}'(z)}=
-\sum_{n\geq1}\frac{2z}{\omega_{a,\nu,n}^2-z^2},$$
and consequently we have
$$\sum_{n\geq 1}\frac{1}{\omega_{a,\nu,n}^2-1}=-\frac{1}{2}\cdot \frac{(2\nu^2-a\nu-1)J_{\nu}(1)+(a-2\nu)J_{\nu}'(1)}{(a-\nu)J_{\nu}(1)+J_{\nu}'(1)}.$$
Now, applying again the fact that \cite[p. 196]{landau} for $\nu>-1$ and $\nu+\gamma\geq0$ the function $\nu\mapsto \gamma_{\nu,n}$ is strictly increasing on $(-1,\infty)$ for $n\in\mathbb{N}$ fixed (where $\gamma_{\nu,n}$ is the $n$th positive root of the equation $\gamma{J}_\nu(z)+zJ_\nu'(z)=0$) we obtain that the function $$\nu\mapsto 1-\sum_{n\geq 1}\frac{1}{\omega_{a,\nu,n}^2-1}$$ is strictly increasing on $(-1,\infty).$
Consequently we have that $$\sum_{n\geq 1}\frac{1}{\omega_{a,\nu,n}^2-1}\leq 1$$
if and only if $\nu\geq \nu_a,$ where $\nu_a$ is the unique root of the equation
$$(a-2\nu+2)J_{\nu}'(1)+(2\nu^2-a\nu-2\nu+2a-1)J_{\nu}(1)=0,$$
or equivalently
$$(2a-1)J_{\nu}(1)-(a-2\nu+2)J_{\nu+1}(1)=0,$$
on $(-1,\infty).$ On the other hand, we know that $\omega_{a,\nu,1}^2>\frac{4a(\nu+1)}{a+2}$ for $\nu>-1$ and $a>0,$ see \cite[Theorem 6.1]{ismail}. This in turn implies that $\omega_{a,\nu,1}>1$ for $\nu>-\frac{3}{4}$ and $a\geq\frac{2}{4\nu+3}.$ Thus, by using Lemma \ref{lem1} we conclude that the function $w_{a,\nu}$ is indeed starlike and all of its derivatives are close-to-convex in $\mathbb{D}$ if and only if $\nu\geq\nu_a.$
\end{proof}

\end{document}